\documentclass[conference] {IEEEtran}

\usepackage{enumerate,xypic,mathdots,mathrsfs,amssymb,amsmath,cases}

\newtheorem{prop}{Proposition}
\newtheorem{theorem}{Theorem}
\newtheorem{cor}{Corollary}
\newtheorem{lemma}{Lemma}
\newtheorem{defn}{Definition}

\newcommand{\A}{\mathcal{A}}

\newcommand{\Norm}[1]{\left\Vert #1 \right\Vert}

\renewcommand{\subset}{\subseteq}

\DeclareMathOperator{\tr}{tr}

\DeclareMathOperator{\spn}{span}
\DeclareMathOperator{\ind}{ind}

\begin{document}

\title{The Zak transform and representations induced from characters of an abelian subgroup}
\author{
\IEEEauthorblockN{Joseph W.\ Iverson}
\IEEEauthorblockA{Department of Mathematics \\ Iowa State University \\ Ames, IA 50011 USA  \\ Email: jwi@iastate.edu}
}
\date{\today}
\maketitle

\pagestyle{plain}

\begin{abstract}
We consider a variant of the Zak transform for a finite group $G$ with respect to a fixed abelian subgroup $H$, and demonstrate a relationship with representations of $G$ induced from characters of $H$.
We also show how the Zak transform can  be used to study right translations by $H$ in $L^2(G)$, and give some examples of applications for equiangular tight frames.
\end{abstract}

The purpose of this note is to demonstrate some connections between the Zak transform and the theory of induced representations.
We also show some applications for right shift-invariant spaces, and for equiangular tight frames that occur as orbits of induced representations.
For the sake of clarity, we restrict our attention to the setting of finite groups, where we can safely ignore convergence issues. 
However, the results in Sections~\ref{sec:I} and~\ref{sec:II} should also hold (with suitable modification) on locally compact groups.

\section{Induced representations}
\label{sec:I}

Fix a finite group $G$, an abelian subgroup $H$, and a transversal $\Omega \subset G$ for $G/H$.
Denote $\hat{H}$ for the Pontryagin dual group of characters $\alpha \colon H \to \mathbb{T}$ under pointwise multiplication. 
We equip each of $G,H,\Omega$ with counting measure, while $\hat{H}$ is given probability measure. 
The left regular representation of $G$ on $L^2(G)$ is given by $(L_x f)(y) = f(x^{-1} y)$, where $x, y \in G$ and $f \in L^2(G)$.
Given a Hilbert space $\mathcal{H}$, we write $L^2(\hat{H};\mathcal{H})$ for the space of functions $\varphi \colon \hat{H} \to \mathcal{H}$, with the inner product
\[ \langle \varphi,\psi \rangle_{L^2(\hat{H}; \mathcal{H})} := \frac{1}{|\hat{H}|} \sum_{\alpha \in \hat{H}} \langle \varphi(\alpha), \psi(\alpha) \rangle_{\mathcal{H}}. \] 

With each $\alpha \in \hat{H}$ we associate the Hilbert space
\[ \mathcal{F}_\alpha = \{ f \in L^2(G) : f(xh) = \overline{\alpha(h)} f(x) \ \forall x\in G,\, h \in H\}, \]
\[ \langle f, g \rangle_{\mathcal{F}_\alpha} := \frac{1}{|H|} \sum_{x \in G} f(x) \overline{g(x)} = \sum_{x \in \Omega} f(x) \overline{g(x)}. \]
Notice that any $f \in \mathcal{F}_\alpha$ is completely determined by its values on $\Omega$, and restriction defines a unitary $\mathcal{F}_\alpha \cong L^2(\Omega)$.
The \emph{induced representation} $\ind_H^G \alpha \colon G \to U(\mathcal{F}_\alpha)$ is given by
\[ [ (\ind_H^G \alpha) (x) f](y) := f(x^{-1} y) \qquad (f \in \mathcal{F}_\alpha,\ x,y \in G). \]
Intuitively, $\ind_H^G \alpha$ performs time-frequency shifts on $L^2(\Omega)$: with respect to the standard basis, each $(\ind_H^G \alpha)(x)$ is given by a monomial (``phased permutation'') matrix.

We will consider the induced representations as an ensemble akin to $\bigoplus_{\alpha \in \hat{H}} \ind_H^G \alpha$.
To that end, the direct integral $\int_{\hat{\alpha}}^{\oplus} \mathcal{F}_\alpha\, d\alpha$ is the Hilbert space consisting of functions ${\varphi \colon \hat{H} \to L^2(G)}$ such that $\varphi(\alpha) \in \mathcal{F}_\alpha$ for every $\alpha \in \hat{H}$, equipped with the inner product
\[ \langle \varphi,\psi \rangle_{\int_{\hat{\alpha}}^{\oplus} \mathcal{F}_\alpha\, d\alpha} := \frac{1}{|\hat{H}|} \sum_{\alpha \in \hat{H}} \langle \varphi(\alpha), \psi(\alpha) \rangle_{\mathcal{F}_\alpha}. \]
Overall, $\int_{\hat{H}}^\oplus \mathcal{F}_\alpha\, d\alpha$ differs from $\bigoplus_{\alpha \in \hat{H}} \mathcal{F}_\alpha$ only by a factor of $1/|\hat{H}|$ in its inner product.
As with $\bigoplus_{\alpha \in \hat{H}} \ind_H^G \alpha$,
the representation $\int_{\hat{H}}^\oplus \ind_H^G \alpha \, d\alpha \colon G \to U(\int_{\hat{H}}^\oplus \mathcal{F}_\alpha \, d\alpha)$ is defined by applying $\ind_H^G \alpha$  coordinate-wise: for $x \in G$ and $\varphi \in \int_{\hat{H}}^\oplus \mathcal{F}_\alpha\, d\alpha$, $( \int_{\hat{H}}^\oplus \ind_H^G \alpha \, d\alpha) (x) \varphi \in \int_{\hat{H}}^\oplus \mathcal{F}_\alpha\, d\alpha$ is the vector-valued function on $\hat{H}$ whose value at $\beta $ is
\[ [( \int_{\hat{H}}^\oplus \ind_H^G \alpha \, d\alpha) (x) \varphi](\beta) := (\ind_H^G\beta) (x)\, [\varphi(\beta)]. \]

In~\cite{I}, we extended the Zak transform construction of Weil~\cite{W} to analyze the left regular representation of a locally compact group, restricted down to an abelian subgroup.
For the specific subgroup $\mathbb{Z} \leq \mathbb{R}$, the construction of~\cite{I} reduces to what is usually called the Zak transform $L^2(\mathbb{R}) \to L^2([0,1]^2)$, which features prominently in time-frequency analysis~\cite{Gr:2001}.
In the present setting, the Zak transform of~\cite{I} amounts to the operator $Z' \colon L^2(G) \to L^2(\hat{H}; L^2(G))$ given by
\[ (Z' f)(\alpha)(x) := \sum_{h\in H} f(hx) \overline{\alpha(h)}, \]
where $f \in L^2(G)$, $\alpha \in \hat{H}$, and $x \in G$.
Here, we introduce a modified version more suitable for the analysis of right translations by $H$ (as opposed to the left translations in~\cite{I}).

\begin{defn}
Given any $f \in L^2(G)$ and $x \in G$, define $f_x \in L^2(H)$ by $f_x(h) = f(xh)$.
The \emph{Zak transform} of $f$ evaluated at $\alpha \in \hat{H}$ is the function $(Zf)(\alpha) \in L^2(G)$ with
\[ (Zf)(\alpha)(x) := \hat{ f_x}(\alpha^{-1}) = \sum_{h\in H} f(xh) \alpha(h) \qquad (x\in G). \]
\end{defn}

This formula appears widely in the basic theory of induced representations on locally compact groups~\cite{F,KT}.
The point here is that we obtain a version of the Zak transform by allowing $\alpha$ to vary across $\hat{H}$.


\begin{theorem}
\label{thm:ind}
The Zak transform $Z \colon L^2(G) \to \int_{\hat{H}}^\oplus \mathcal{F}_\alpha\, d\alpha$ is a unitary operator that intertwines the left regular representation of $G$ with $\int_{\hat{H}}^\oplus \ind_H^G \alpha \, d\alpha$.
\end{theorem}

\begin{IEEEproof}
It is easy to see that $Z$ maps $L^2(G)$ into $\int_{\hat{H}}^\oplus \mathcal{F}_\alpha\, d\alpha$ while intertwining the left regular representation of $G$ with $\int_{\hat{H}}^\oplus \ind_H^G \alpha\, d\alpha$.
A straightforward application of Plancherel's Theorem on $L^2(H)$ shows that $Z$ is an isometry.
By counting dimensions, we conclude that $Z$ is unitary.
\end{IEEEproof}

\section{Shift-invariant spaces}
\label{sec:II}

How about the right regular representation?
We now show the Zak transform diagonalizes right translations by $H$.
To that end, it will be convenient to identify $\mathcal{F}_\alpha \cong L^2(\Omega)$ by restriction.
(Recall $\Omega$ is a transversal for $G/H$.)
We can then view $Z$ as a function $Z_r \colon L^2(G) \to L^2(\hat{H};L^2(\Omega))$ by restricting $(Z_r f)(\alpha) := \left. (Zf)(\alpha) \right|_{\Omega}$.

In the theorem below, we write $R_h \in U(L^2(G))$ for right translation by $h \in H$, namely, $(R_h f)(x) := f(xh)$. Meanwhile, the modulation representation of $H$ on $L^2(\hat{H};L^2(\Omega))$ is defined by the formula $(M_h f)(\alpha) := \overline{\alpha(h)}\cdot f(\alpha) $.

\begin{theorem}
\label{thm:mod}
The Zak transform defines a unitary $Z_r \colon L^2(G) \to L^2(\hat{H};L^2(\Omega))$ that intertwines right $H$-translations with modulations: $Z_r(R_h f) = M_h(Z_r f)$.
\end{theorem}

\begin{IEEEproof}
Restriction defines unitaries $\mathcal{F}_\alpha \cong L^2(\Omega)$ for each $\alpha \in \hat{H}$.
Overall, it gives an identification $\int_{\hat{H}}^\oplus \mathcal{F}_\alpha\, d\alpha \cong L^2(\hat{H};L^2(\Omega))$.
That $Z_r$ is unitary is thus a corollary of Theorem~\ref{thm:ind}.
The intertwining formula is obtained by reindexing the sum in the definition of $Z$.
\end{IEEEproof}

\begin{defn}
A subspace $V \subset L^2(G)$ is called \emph{right $H$-shift invariant} (right $H$-SI) if $R_h V = V$ for every $h \in H$.
A \emph{range function} is a mapping $J \colon \hat{H} \to \{ \text{subspaces of }L^2(\Omega)\}$.
\end{defn}

Following the well-trodden path of~\cite{He,S,BDR,B}, the Zak transform serves to characterize right $H$-SI spaces in terms of range functions.
Similar connections between the Zak transform and (left) shift-invariant spaces were first observed in~\cite{I}, \cite{I2}, and, independently,~\cite{BHP2}.

\begin{cor}
\label{cor:SI}
Given a range function $J$, put
\[ V_J := \{ f \in L^2(G) : (Z_rf)(\alpha) \in J(\alpha) \ \forall \alpha \in \hat{H} \}. \]
Then the mapping $J \mapsto V_J$ is a one-to-one correspondence between range functions and right $H$-SI subspaces of $L^2(G)$.
\end{cor}

\begin{IEEEproof}
With each range function $J$, associate the space
\begin{equation}
\label{eq:MJ}
\mathcal{M}_J = \{ f \in L^2(\hat{H};L^2(\Omega)) : f(\alpha) \in J(\alpha) \ \forall \alpha \in \hat{H} \}.
\end{equation}
Thus, $Z_rV_J = \mathcal{M}_J$.
In the language of \cite{BR}, the characters of $\hat{H}$ are a ``determining set'' for $L^1(\hat{H})$, by~\cite[Lemma~5.2]{I}.
Applying~\cite[Theorem~2.4]{BR}, we see that $J \mapsto \mathcal{M}_J = Z_rV_J$ gives a one-to-one correspondence between range functions and modulation-invariant (MI) spaces in $L^2(\hat{H};L^2(\Omega))$.
By Theorem~\ref{thm:mod}, $Z_r^{-1}$ gives a one-to-one correspondence between MI spaces in $L^2(\hat{H};L^2(\Omega))$ and right $H$-SI spaces in $L^2(G)$.
Composing these correspondences proves the corollary.
\end{IEEEproof}

\begin{defn}
In a Hilbert space $\mathcal{H}$, a sequence of vectors $\Phi = \{ f_j \}_{j \in I}$ is a \emph{frame} if there are \emph{bounds} $A,B > 0$ such that $A \Norm{g}^2 \leq \sum_{j \in I} | \langle g, f_j \rangle |^2 \leq B \Norm{g}^2$ for every $g \in \mathcal{H}$.
It is a \emph{tight frame} if $A=B$, and \emph{equiangular} if there is a constant $C$ such that $| \langle f_i, f_j \rangle | = 1$ when $i=j$ and $C$ otherwise. An \emph{equiangular tight frame} (ETF) is both tight and equiangular.
\end{defn}

Given a finite sequence $\mathscr{A} = \{ f_j \}_{j\in [n]}$ in $L^2(G)$, we denote $E(\A) = \{ R_h f_j \}_{h\in H,\, j \in [n]}$ for the sequence of its right shift by $H$, and $S(\A) = \spn E(\A)$ for the right $H$-SI space it generates.
Applying Corollary~\ref{cor:SI}, it is easy to see that $S(\A) = V_{J_{\mathscr{A}}}$, where $J_{\mathscr{A}}(\alpha) := \spn\{ (Z_r f_j)(\alpha) : j \in [n] \}$.

\begin{cor}
\label{cor:frm}
Let $\mathscr{A} = \{ f_j \}_{j\in [n]}$ be a finite sequence in $L^2(G)$. For any $A,B > 0$, the following are equivalent:
\begin{itemize}
\item[(i)]
$E(\A)$ is a frame for $S(\A)$ with bounds $A,B$.

\item[(ii)]
For every $\alpha \in \hat{H}$, $\{ (Z_rf_j)(\alpha) \}_{j\in [n]}$ is a frame for $J_{\mathscr{A}}(\alpha)$ with bounds $A,B$.
\end{itemize}
\end{cor}

Similar results appear in~\cite{BL:93,B,RS,I,I2,BHP2}.

\begin{IEEEproof}
Let $\mathcal{M}_{J_{\mathscr{A}}}$ be as in~\eqref{eq:MJ}.
By Theorem~\ref{thm:mod}, (i) is equivalent to
\begin{itemize}
\item[(i')]
$\{ M_h (Z_rf_j) \}_{h\in H,\, j \in [n]}$ is a frame for $M_{J_\mathscr{A}}$ with bounds $A,B$.
\end{itemize}
For each $h \in H$, let $e_h \in L^\infty(\hat{H})$ be the evaluation character given by $e_h(\alpha) = \overline{\alpha(h)}$.
In the language of~\cite{I}, $\mathcal{D}:= \{ e_h \}_{h \in H}$ is a ``Parseval determining set'' for $L^1(\hat{H})$, by~\cite[Lemma~5.2]{I}.
Then \cite[Theorem 2.10]{I} gives the equivalence of (i') and (ii).
\end{IEEEproof}

In the case of a single generator $\mathscr{A} = \{f\}$, Corollary~\ref{cor:frm} reduces to the following analogue of \cite[Theorem 3]{BL:93}.
\begin{cor}
For any $f \in L^2(G)$ and any $A,B > 0$, $\{ R_h f \}_{h\in H}$ is a frame for $S(\{f\})$ with bounds $A,B$ if and only if the following holds for every $\alpha \in \hat{H}$:
\[ \sum_{x \in \Omega} | (Z_r f)(\alpha)(x) |^2 \in \{0\} \cup [A,B]. \]
\end{cor}

\section{Equiangular tight frames}

For the remainder of the paper, we focus on the special case where $G = H \rtimes K$.
Then $G$ acts on $\hat{H}$ by the formula $(x\cdot \alpha)(h) = \alpha(x^{-1} h x)$.
We choose $\Omega = K$ for a transversal of $G/H$.
Given $\alpha \in \hat{H}$, we can identify $\mathcal{F}_\alpha \cong L^2(K)$ by restriction, as above.
In that case, $\ind_H^G \alpha$ gives time-frequency shifts in $L^2(K)$: for $h \in H$, $x,y \in K$, and $f \in L^2(K)$,
\[ [(\ind_H^G \alpha)(hx) f](y) = (y \cdot \alpha)(h)\, f(x^{-1} y).  \]
By the Mackey machine, $\ind_H^G \alpha$ is irreducible whenever the \emph{little group} $K_\alpha := \{ x \in K : x\cdot \alpha = \alpha \}$ is trivial (see Section~6.6 of~\cite{F}).

In this section, we study ETFs for $L^2(K)$ whose vectors span lines that occur as orbits of $\ind_H^G \alpha$.
To that end, it will be helpful to define the \emph{projective stabilizer} of a unit vector $f \in L^2(K)$ to be
the group of all $x \in G$ stabilizing $\spn\{f\}$,
\[ L = \{ x \in G : (\ind_H^G \alpha)(x)f = c f \text{ for some }c \in \mathbb{T} \}. \]
The \emph{function of positive type} associated with $f \in L^2(K)$ is $g\in L^2(G)$ given by $g(x) = \langle f, (\ind_H^G \alpha)(x) f \rangle$.

\begin{defn}
Given a finite sequence $\Phi = \{ f_j \}_{j \in [n]}$ of unit vectors, a \emph{projective reduction} of $\Phi$ is a subsequence $\Phi' $ such that every line in the set $\{ \spn\{f_j\}: j \in [n] \}$ is represented exactly once in $\Phi'$.
\end{defn}

Projective reduction is not unique, but any two projective reductions of $\Phi$ are equiangular (resp.\ tight) at the same time.

\begin{lemma}
\label{lem:pos}
Assume $G = H \rtimes K$, and fix $\alpha \in \hat{H}$.
Given a unit vector $f \in L^2(K)$,  let $g \in L^2(G)$ be its function of positive type. 
Then any projective reduction $\Phi'$ of $\Phi := \{ (\ind_H^G \alpha)(x) f \}_{x \in G}$ is equiangular if and only if there exists $C \neq 1$ and a set $L\subset G$ containing $1_G$ such that
\begin{equation}
\label{eq:EqPos}
|g(x)|^2 = \begin{cases}
1, & \text{if } x \in L; \\
C, & \text{otherwise.}
\end{cases}
\end{equation}
In that case, $L$ is the projective stabilizer of $f$, and $\Phi'$ contains $|G|/|L|$ vectors.
Moreover, $\Phi'$ is an ETF for $L^2(K)$ if and only if $f$ is cyclic for $\ind_H^G \alpha$  and \eqref{eq:EqPos} holds with $C = \frac{|H| - |L|}{|G| - |L|}$.
\end{lemma}

\begin{IEEEproof}
We abbreviate $\pi = \ind_H^G \alpha$ throughout.
By the equality condition in Cauchy-Schwarz, the projective stabilizer of $f$ is the set $L$ of all $x \in G$ such that $|g(x)|^2 = 1$.
Fix a transversal $\{ x_j \}_{j\in [n]}$ for $G/L$.
For any $x,y \in G$, we have $\spn\{\pi(x) f\} = \spn\{\pi(y) f\}$ if and only if there exists $c\in \mathbb{T}$ such that $\pi(y^{-1} x) f = cf$, if and only if $xL = yL$.
Thus, $\Phi' := \{ \pi(x_j) f \}_{j \in [n]}$ is a projective reduction of $\Phi$.

For any $j \in [n]$ and any $y \in L$ with $\pi(y)f = cf$, we have 
\[ |g(x_j y)|^2 = | \langle f, \pi(x_j) c f \rangle |^2 = |g(x_j)|^2, \]
so $|g|^2$ is a function on $G/L$.
Moreover, for any $i \neq j$, $|\langle \pi(x_i) f, \pi(x_j) f \rangle |^2 = | g(x_i^{-1} x_j)|^2$.
Therefore, 
\[ \{ | \langle \pi(x_i) f, \pi(x_j) f \rangle |^2 : i \neq j \} = \{ |g(x)|^2 : x \notin L \}. \]
It follows that $\Phi'$ is equiangular if and only if~\eqref{eq:EqPos} holds.

Finally, assume $\Phi'$ is equiangular, and put $d = |K|$.
Then $\Phi'$ is a frame for $L^2(K)$ if and only if $L^2(K) = \spn \Phi' = \spn \Phi$, i.e. $f$ is a cyclic vector.
In that case, the Welch bound~\cite{Wel74} states that $\Phi'$ is an ETF for $L^2(K)$ if and only if 
\[ | \langle \pi(x_i) f, \pi(x_j) f \rangle |^2 = C = \frac{n-d}{d(n-1)} = \frac{ |H| - |L| }{ |G| - |L| } \]
whenever $i \neq j$.
\end{IEEEproof}

In many cases of interest, it is simpler to understand $|g|^2$ by applying the Zak transform with the following formulae.

\begin{lemma}
\label{lem:ZtPos}
Assume $G = H \rtimes K$, and fix $\alpha \in \hat{G}$.
Given $f \in L^2(K)$, let $g \in L^2(G)$ be its function of positive type.
Then for any $\beta \in \hat{H}$ and $x \in K$,
\begin{equation}
\label{eq:ZtPos}
(Z_rg)(\beta)(x) = |H| \sum_{\substack{y \in K, \\ y\cdot \alpha = \beta}} f(xy) \overline{f(y)}.
\end{equation}
\end{lemma}

\begin{IEEEproof}
It follows easily from the definitions by applying character orthogonality on $H$.
\end{IEEEproof}

\begin{lemma}
\label{lem:conv}
For any $f, g \in L^2(G)$ and any $\alpha \in \hat{H}$, we have
\[ [Z(f\cdot g)](\alpha) = \frac{1}{|\hat{H}|} \sum_{\beta \in \hat{H}} (Zf)(\beta) \cdot (Zg)(\beta^{-1} \alpha) \]
and $(Z \overline{f})(\alpha) = \overline{ (Zf)(\alpha^{-1}) }$.
In particular,
\begin{equation}
\label{eq:Zf2}
(Z |f|^2)(\alpha) = \frac{1}{|\hat{H}|} \sum_{\beta \in \hat{H}} (Zf)(\beta) \cdot \overline{ (Zf)(\alpha^{-1} \beta) }.
\end{equation}
(Here, multiplication and complex conjugation in $L^2(G)$ are interpreted pointwise.)
\end{lemma}

\begin{IEEEproof}
It follows easily from the fact that the Fourier transform $L^2(H) \to L^2(\hat{H})$ intertwines pointwise multiplication and conjugation with convolution and involution, respectively. 
We leave details to the reader.
\end{IEEEproof}

\subsection{Affine linear groups and Paley difference sets}

Fix a prime power $q = p^r > 3$ with $q \equiv 3 \mod 4$, and let $\mathbb{F}_q$ be the finite field of order $q$.
Write $\mathbb{F}_q^{\times 2}$ for the group of nonzero quadratic residues. 
Take $G$ to be the group of all affine transformations $x \mapsto ax+b$ on $\mathbb{F}_q$ with $a \in \mathbb{F}_q^{\times 2}$ and $b \in \mathbb{F}_q$.
Then $G = H \rtimes K$, where $H$ is the group of translations $\tau_b(x) := x+b$ and $K$ that of dilations $\theta_a(x) := ax$.
%
If $\omega = \exp(2\pi i/p)$, then the dual group $\hat{H}$ consists of all characters $\alpha_b(\tau_c) := \omega^{\tr(bc)}$, where $\tr \colon \mathbb{F}_q \to \mathbb{F}_p$ is the field trace.
(See~\cite{LN83}.)
Then $\alpha_a \alpha_b = \alpha_{a+b}$, and the action of $K$ on $\hat{H}$ is given by $\vartheta_a \cdot \alpha_b = \alpha_{a b}$.

\begin{prop}
Let $f \in L^2(K)$ be the constant function $f(\vartheta_a) \equiv \sqrt{2/(q-1)}$.
Then the projective reduction of $\Phi = \{ (\ind_H^G \alpha_1)(x) g \}_{x\in G}$ is an ETF of $q$ vectors in $L^2(K)$, a space of dimension $(q-1)/2$.
\end{prop}

The resulting ETF is well known, but not by this construction.
In fact, direct examination of the short fat matrix representing $\Phi'$ shows we obtain the harmonic ETF~\cite{SH,XZG,DF} corresponding to the Paley difference set $\mathbb{F}_q^{\times 2}$ in $\mathbb{F}_q$.

\begin{IEEEproof}
Let ${g \in L^2(G)}$ be the function of positive type associated with $f$.
For any $\beta = \alpha_b \in \hat{H}$ and any $x = \vartheta_a \in K$, we have $x\cdot \alpha_1 = \beta$ if and only if $a = b$.
Comparing with~\eqref{eq:ZtPos}, we deduce that $(Z_r g)(\alpha_b) \equiv  \frac{2q}{q-1}\, \chi_{\mathbb{F}_q^{\times 2}}(b)$.
Then \eqref{eq:Zf2} produces
\[ (Z_r |g|^2)(\alpha_b) \equiv \frac{4q}{(q-1)^2} \sum_{c \in \mathbb{F}_q} \chi_{\mathbb{F}_q^{\times 2}}(c)\, \chi_{\mathbb{F}_q^{\times 2}}(c-b). \]
It is well known that $\mathbb{F}_q^{\times 2}$ is a $(q,\frac{q-1}{2},\frac{q-3}{4})$-difference set in $(\mathbb{F}_q, +)$, so that
\[ \sum_{c \in \mathbb{F}_q} \chi_{\mathbb{F}_q^{\times 2}}(c)\, \chi_{\mathbb{F}_q^{\times 2}}(c-b) = 
\begin{cases}
\frac{q-1}{2}, & \text{if } b = 0; \\[3 pt]
\frac{q-3}{4}, & \text{otherwise.}
\end{cases} \]
(See~\cite{CRC} for background.) Overall, $(Z_r |g|^2)(\alpha_b) \equiv 2q/(q-1)$ when $b=0$, and $q(q-3)/(q-1)^2$ otherwise.

By comparison, let $g' \in L^2(G)$ be the desired value of $|g|^2$,
\[ g'(x) = \frac{q+1}{(q-1)^2} + \frac{q(q-3)}{(q-1)^2}\chi_K = \begin{cases} 1, & \text{if }x \in K; \\ \frac{q+1}{(q-1)^2}, & \text{otherwise.} \end{cases} \]
The constant function $1 \in L^2(G)$ has Zak transform $(Z_r 1)(\alpha_b)(x) = \sum_{h \in H} \alpha_b(h) \equiv q\delta_{0,b}$, while
\[ (Z_r \chi_K)(\alpha_b)(x) = \sum_{h \in H} \chi_K(hx) \alpha_b(h) \equiv 1. \]
By linearity, $(Z_r g')(\alpha_b) \equiv \frac{q(q+1)}{(q-1)^2} \delta_{0,b} + \frac{q(q-3)}{(q-1)^2} = 2q/(q-1)$ when $b=0$, and $q(q-3)/(q-1)^2$ otherwise. Since $Z_r$ is injective, $|g|^2 = g'$.

Finally, the little group $K_{\alpha_1}$ is trivial, so $\ind_H^G \alpha_1$ is irreducible, and $f \neq 0$ is a cyclic vector.
By Lemma~\ref{lem:pos}, the projective reduction of $\Phi = \{ (\ind_H^G \alpha_1)(x) f \}_{x\in G}$ is an ETF
consisting of $q$ vectors.
\end{IEEEproof}

\subsection{Finite Heisenberg groups and SIC-POVMs}

Fix an integer $d \geq 2$, and define
\[ G = \langle r, s, t : r^d = s^d = t^d = [r,s] = [t,s] = 1, trt^{-1} = rs \rangle. \]
This is the Heisenberg group {mod $d$}.
We have $G = H \rtimes K$, where $H := \langle r,s \rangle \cong \mathbb{Z}_d \times \mathbb{Z}_d$ and $K := \langle t \rangle \cong \mathbb{Z}_d$.
Denoting $\omega = e^{2\pi i /d}$, the dual group $\hat{H}$ consists of characters $\alpha_{a,b}(r^m s^n) := \omega^{am+bn}$, with $a,b \in \mathbb{Z}_d$.
In this notation, $\alpha_{a,b} \alpha_{a',b'} = \alpha_{a+a',b+b'}$, and the action of $K$ on $\hat{H}$ satisfies $t^k \cdot \alpha_{a,b} = \alpha_{a-kb,b}$.

Define $\pi := \ind_H^G \alpha_{0,1}$.
Then the little group $K_{\alpha_{0,1}}$ is trivial, and $\pi$ is irreducible.
The projective stabilizer of every $f \in L^2(K)$ contains $L := \langle s \rangle$, since $\pi(s^n) = \omega^n I$.
If $L$ is the entire projective stabilizer of $f$, then the projective reduction of $ \{ \pi(x) f \}_{x\in G}$ contains $d^2$ vectors in a space of dimension $d$.
Any ETF of this form is known as a \emph{symmetric informationally complete positive operator-valued measure} (SIC-POVM) in quantum information theory~\cite{RBS}.
\emph{Zauner's conjecture} posits that SIC-POVMs exist for every $d$~\cite{Zauner:99}.
A large body of numerical evidence supports this conjecture~\cite{SG:2010}.

The following characterization of SIC-POVMs generated by $\pi$ has been found many times~\cite{ADF,BW,Kh08}.
We give a simple proof using the Zak transform.

\begin{prop}
Let $f \in L^2(K)$ be an arbitrary unit vector.
Then the projective reduction of $\Phi := \{ \pi(x) f \}_{x\in G}$ is a SIC-POVM if and only if the following holds for every $a,h \in \mathbb{Z}_d$:
\[ 
\sum_{b \in \mathbb{Z}_d} f(t^{h+b}) \overline{f(t^b)} \overline{f(t^{h+a+b})} f(t^{a+b}) = \frac{\delta_{a,0} + \delta_{h,0}}{d+1}.
\]
\end{prop}

\begin{IEEEproof}
Let $g \in L^2(G)$ be the function of positive type associated with $f$.
Take any $a,c,h \in \mathbb{Z}_d$.
Since $t^k \cdot \alpha_{0,1} = \alpha_{-k,1}$, Lemma~\ref{lem:ZtPos} says that $(Z_r g)(\alpha_{a,c})(t^h) = d^2 f(t^{h-a}) \overline{ f(t^{-a}) } \delta_{c,1}$.
Applying Lemma~\ref{lem:conv}, we obtain
\begin{equation}
\label{eq:SIC1}
(Z_r |g|^2)(\alpha_{a,c})(t^h) \qquad \qquad \qquad \qquad \qquad \qquad \qquad \qquad
\end{equation}
\[ = d^2 \sum_{b \in \mathbb{Z}_d} f(t^{h-b}) \overline{f(t^{-b})} \overline{f(t^{h-b+a})} f(t^{a-b}) \delta_{c,0}. \]

On the other hand, let $g' = \frac{1}{d+1} + \frac{d}{d+1} \chi_L \in L^2(K)$, so that
$g'(x) = 1$ when $x \in L$, and $1/(d+1)$ otherwise.
By Lemma~\ref{lem:pos}, the projective reduction of $\Phi$ is a SIC-POVM if and only if $|g|^2 = g'$.
As in the previous example, the constant function $1 \in L^2(G)$ has Zak transform $(Z_r 1)(\alpha_{a,c})(t^h) \equiv d^2 \delta_{a,0} \delta_{c,0}$.
Meanwhile,
\begin{align*}
(Z_r \chi_L)(\alpha_{a,c})(t^h) &= \sum_{m,n \in \mathbb{Z}_d} \chi_L(t^h r^m s^n) \alpha_{a,c}(r^m s^n) \\
&= \sum_{m,n \in \mathbb{Z}_d} \chi_L(r^m s^{hm+n} t^h) \alpha_{a,c}(r^m s^n) \\
&= \delta_{h,0} \sum_{n \in \mathbb{Z}_d} \alpha_{a,c}(s^n) \\
&= d \delta_{h,0} \delta_{c,0}.
\end{align*}
By linearity,
\begin{equation}
\label{eq:SIC2}
(Z_r g')(\alpha_{a,c})(t^h) = \frac{d^2}{d+1} \delta_{c,0}(\delta_{a,0} + \delta_{h,0}).
\end{equation}
The proposition follows by comparing \eqref{eq:SIC1} and \eqref{eq:SIC2}.
\end{IEEEproof}


\bibliographystyle{IEEEtran}
\bibliography{IEEEabrv,SAMPTA2019}

\end{document}